\documentclass[graybox]{svmult}

\usepackage{helvet}         
\usepackage{courier}        

\usepackage{graphicx}              

\usepackage{siunitx}
\usepackage{amsbsy}
\usepackage{amsmath,bm}
\usepackage{amssymb}
\usepackage{amsfonts}

\usepackage{multirow}

\usepackage[bottom]{footmisc}
\usepackage{cite}

\usepackage{url}

\def\d{\mbox{d}}
\def\dt{\mbox{d}t}
\def\dtau{\mbox{d}\tau}
\newcommand{\ddt}{\frac{\d}{\dt}}
\newcommand{\interval}{\mathcal{I}}

\if 11
\newcommand\myopbf\bf                     
\newcommand\myboldmath\boldsymbol
\newcommand{\vect}[1]{{\mathbf {#1}}} 		
\else
\newcommand\myopbf\rm                           
\newcommand{\vect}[1]{#1}		
\newcommand\myboldmath\unboldmath
\fi

\if 0 
\newcommand{\Atop}[1][{}]{\vect{A}_{\!#1}}
\newcommand{\jcurrv}[1][{}]{\myboldmath\jmath_{#1}}  %
\newcommand{\npot}[1][{}]{{\vect{e}}_{{#1}}}
\newcommand{\charge}[1][{}]{{\vect{q}}_{{#1}}}

\newcommand{\flux}[1][{}]{\myboldmath\phi_{{#1}}}

\newcommand{\vsrc}[1][{}]{{\vect{v}}_{{#1}}}
\newcommand{\csrc}[1][{}]{{\myboldmath\imath}_{{#1}}}
\newcommand{\resi}[1][{}]{{\vect{g}}}

\else
\newcommand{\Atop}[1][{}]{{A}_{#1}}
\newcommand{\jcurrv}[1][{}]{\jmath_{#1}}

\newcommand{\npot}[1][{}]{{{e}}_{{#1}}}
\newcommand{\charge}[1][{}]{{q}_{{#1}}}

\newcommand{\flux}[1][{}]{\phi_{{#1}}}

\newcommand{\vsrc}[1][{}]{{{v}}_{{#1}}}
\newcommand{\csrc}[1][{}]{{\imath}_{{#1}}}
\newcommand{\resi}[1][{}]{{{g}}}

\renewcommand{\vect}[1]{#1}		
\renewcommand\myboldmath\unboldmath
\fi

\newcommand{\tim}{t}

\newcommand{\tr}{\top}

\usepackage{bbm}
\newcommand{\real}{\mathbbm{R}}
\newcommand{\diag}{\text{diag}}

\begin{document}

\title*{Port-Hamiltonian Systems Modelling  in Electrical Engineering}
\author{A. Bartel 
\and
M. Clemens 
\and
M. Günther \inst{}
\and
B. Jacob
\and
T. Reis
}
\institute{Andreas Bartel, Michael Günther  \at Bergische Universität Wuppertal, IMACM, Chair of Applied Mathematics, Gaußstraße 20, D-42119 Wuppertal, 
\email{[bartel,guenther]@uni-wuppertal.de}
\and Markus Clemens \at Bergische Universität Wuppertal, IMACM, Chair of Electromagnetic Theory, Rainer-Gruenter-Straße 21, D-42119 Wuppertal, 
  \email{clemens@uni-wuppertal.de}
  \and Birgit Jacob  \at Bergische Universität Wuppertal, IMACM, Chair of Functional Analysis, Gaußstraße 20, D-42119 Wuppertal, 
\email{bjacob@uni-wuppertal.de}
  \and Timo Reis \at TU Ilmenau, Fakultät für Mathematik und Naturwissenschaften, Chair of System Theory and PDEs, PF 10 05 65, D-98684 Ilmenau, \email{timo.reis@tu-ilmenau.de}}

%
%
\maketitle

\abstract*{The port-Hamiltonian modelling framework allows for models that preserve essential physical properties such as energy conservation or dissipative inequalities. If all subsystems are modelled as port-Hamiltonian systems and the inputs are related to the output in a linear manner, the overall system can be modelled as a port-Hamiltonian system (PHS), too, which preserves the properties of the underlying subsystems. If the coupling is given by a skew-symmetric matrix, as usual in many applications, the overall system can be easily derived from the subsystems without the need of introducing dummy variables and therefore artificially increasing the complexity of the system.
Hence the PHS framework is especially suitable for modelling multiphysical systems. 
\\
In this paper we show that port-Hamiltonian systems are a natural generalization of Hamiltonian systems, define coupled port-Hamiltonian systems as ordinary and differential-algebraic equations. To highlight the suitability for electrical engineering applications we derive PHS models for MNA network equations, electromagnetic devices and coupled systems thereof. }

\abstract{The port-Hamiltonian modelling framework allows for models that preserve essential physical properties such as energy conservation or dissipative inequalities. If all subsystems are modelled as port-Hamiltonian systems and the inputs are related to the output in a linear manner, the overall system can be modelled as a port-Hamiltonian system (PHS), too, which preserves the properties of the underlying subsystems. If the coupling is given by a skew-symmetric matrix, as usual in many applications, the overall system can be easily derived from the subsystems without the need of introducing dummy variables and therefore artificially increasing the complexity of the system.
Hence the PHS framework is especially suitable for modelling multiphysical systems. 
\\
In this paper, we show that port-Hamiltonian systems are a natural generalization of Hamiltonian systems, define coupled port-Hamiltonian systems as ordinary and differential-algebraic equations. To highlight the suitability for electrical engineering applications, we derive PHS models for MNA network equations, electromagnetic devices and coupled systems thereof. }

\section{Port-Hamiltonian Systems Modelling in a Nutshell}
Port-Hamiltonian Systems (PHS) are a generalization of Hamiltonian systems 
   \begin{align}
   \label{bcgjr_ham} 
         \dot x & = 
         J \cdot \nabla H(x), \quad x(0)=x_0
     \end{align}
     with $x=(p,q)$ consisting of generalized position $q(t) \in \mathbb R^n$ and momentum $p(t) \in \mathbb R^n$ (where $t\in [0,T]$), the skew-symmetric matrix $J$ given by
     \begin{align*}
         J = 
         \left[\begin{array}{r@{\;\;}r}
           0 & -I \\ I & 0
         \end{array}
         \right]
     \end{align*}
     and the Hamiltonian $H(x)=H(p,q)=U(p)+V(q)$
     given as the sum of potential and kinetic energy, which maps $\mathbb R^n \times \mathbb R^n \rightarrow \mathbb R$ and is twice continuously differentiable.  
     The Hamiltonian flow $\varphi(t;x_0))$, i.e., the solution of \eqref{bcgjr_ham} at time point $t$, starting at the initial value $x(0)=x_0$, is characterized by four geometric properties:
    
 \begin{enumerate}
     \item {\em Preservation of the Hamiltonian}: 
     \[
       \ddt H(\varphi(t;x_0))= (\nabla H(\varphi(t;x_0)))^\top J (\nabla H(\varphi(t;x_0))) = 0.
      \] 
     
     \item {\em Time-reversibility}: $$ \rho \circ \varphi(t;x_0) \circ \rho \circ \varphi(t;x_0) = x_0,$$ with $\rho(p,q)=(-p,q)$, which is a direct consequence of the $\rho$-reversibility of the Hamiltonian flow: $\rho \circ J \nabla H(\varphi(t;x_0))) = -J \nabla H(\rho \circ \varphi(t;x_0)))$.
     
     \item {\em Symplectic structure of the Hamiltonian flow}:
     \begin{align*}
         \Psi(t)^\top J^{-1} \Psi(t) = J^{-1} , \qquad 
         \textstyle
         \Psi(t):= \frac{\partial \varphi(t;x_0)}{\partial x_0},
   \end{align*}
    which is a direct consequence of the skew-symmetry of $J$.
    
         \item {\em Volume-preservation}: 
         $$(\det \Psi(t))^2=1,$$ 
         which follows immediately from the symplectic structure 
         in 3.
 \end{enumerate}
 \subparagraph{First generalization step: arbitrary skew-symmetric matrices $J$}
If we replace in~\eqref{bcgjr_ham} $J$ by an arbitrary skew-symmetric matrix, the Hamiltonian is still preserved. As $x$ will loose its characterization  as  generalized positions and momenta of classical mechanics, time-reversibility will generally not hold anymore. However, the symplectic structure of the flow still holds 
in the case of a regular $J$, and volume preservation is still a consequence of the Hamiltonian flow.
 
 \subparagraph{Second generalization step: adding dissipation to the system}
 Allowing the flow to become dissipative, we may generalize~\eqref{bcgjr_ham} to the dissipative Hamiltonian system 
\begin{align}
\label{bcgjr_disham}
          \dot x & = 
        (J-R)  \cdot \nabla H(x), \quad x(0)=x_0 
     \end{align}
     with $R\ge 0$ being symmetric and positive 
     semi-definite. In this case, the flow will neither be symplectic nor volume preserving, and the preservation of the Hamiltonian is replaced by the dissipativity condition
    \begin{align*}
         \ddt H(x(t)) & = (\nabla H(x))^\top \dot x = -(\nabla H(x))^\top R \nabla H(x) \le 0 \\
         \Rightarrow
         H(x(t)) & = H(x_0)- \int_0^t (\nabla H(x(\tau)))^\top R \nabla H(x(\tau)) \, \dtau \le H(x_0).
     \end{align*} 
 
 \subparagraph{Third generalization step: coupling to the environment via inputs and outputs}
 Allowing for inputs and outputs to couple the system to the environment, we end up with linear port-Hamiltonian system characterized by 
 \begin{align*}
 \dot x & = 
        (J-R)  \cdot \nabla H(x) + B u(t), \quad x(0)=x_0, \\
        y & = B^\top \nabla H(x) 
     \end{align*}
     with inputs $u(t) \in \mathbb R^p$, outputs $y(t) \in \mathbb R^p$ and port-matrices $B \in \mathbb R^{n \times p}$. The dissipativity inequality now reads
   \begin{align*}
         \ddt H(x(t)) & = (\nabla H(x))^\top \dot x = -(\nabla H(x))^\top R \nabla H(x) + (\nabla H(x))^\top B u(t)) \\
         & = -(\nabla H(x))^\top R \nabla H(x) + y(t)^\top u(t) \le y(t)^\top u(t) \\
         \Rightarrow
         H(x(t)) & = H(x_0)- \int_0^t (\nabla H(x(\tau)))^\top R \nabla H(x(\tau)) \, \dtau
         + \int_0^t  y(\tau)^\top u(\tau) \dtau \\
         & \le H(x_0) + \int_0^t  y(\tau)^\top u(\tau) \dtau.
     \end{align*}    

\subparagraph{Fourth generalization step: PH-DAE systems}
Linear PHS can be easily generalized to PH-DAE systems given by
     \begin{subequations}
     \begin{align}
     \label{eq.lin.ph-dae}
  \ddt (E x) & = 
        (J-R)  \cdot z(x) + B u(t), \quad x(0)=x_0, \\
        y & = B^\top z(x) 
     \end{align}
     \end{subequations}
     with a possibly singular
     matrix $E\in \mathbb{R}^{n\times n}$  and the nonlinear mapping $z: \mathbb R^n \rightarrow \mathbb R^n$ fulfilling the compatibility condition $E^\top z = \nabla H$. Now the dissipativity condition reads 
  \begin{align*}
         H(x(t)) & = H(x_0)- \int_0^t z(x(\tau))^\top R \nabla z(x(\tau)) \, \dtau + \int_0^t  y(\tau)^\top u(\tau) \, \dtau\\
         & \le H(x_0) + \int_0^t  y(\tau)^\top u(\tau) \, \dtau.
     \end{align*}     
      The key point in port-Hamiltonian modelling is the following: there is an easy way to couple arbitrary many PH-DAE system such that the overall system is still a PH-DAE system, which preserves a dissipativity inequality. 
     
     Let us consider $r$ autonomous PH-DAE systems 
 \begin{subequations}
 \label{bcgjr_coupled.phdaes}
\begin{align}
     \ddt (E_i x_i) & = (J_i-R_i) z_i(x_i) 
     + B_i u_i,\\
     y_i & = B_i^{\top} z_i(x_i) 
 \end{align}
   \end{subequations}
  with $r$ Hamiltonians $H_1, H_2, \ldots, H_r$ and compatibility conditions $E_i^\top z_i = \nabla H_i $. 
  If the inputs and outputs fulfill a 
linear interconnection relation $Mu+Ny= 0$ for 
     the aggregated input $u= (u_1,u_2,\ldots,u_r)$ and output $y= (y_1,y_2,\ldots,y_r)$, 
it has been shown in~\cite{MeMo19} that one can write the 
aggregated system as a joint  PH-DAE system  as 
 \begin{align*}
     \ddt \left( \begin{bmatrix} 
     E & 0 & 0   \\ 0 & 0 & 0\\ 0 & 0 & 0\\ 0 & 0 & 0
     \end{bmatrix} 
     \begin{bmatrix} x \\ {\hat u} \\ {\hat y} \end{bmatrix} \right) & =
      \begin{bmatrix}J-R & B & 0 & 0 \\-B^\top & 0
        & I_m & -M^{\top} \\
       0 & -I_m & 0 & -N^{\top} \\
       0 & M & N & 0\end{bmatrix}
          \begin{bmatrix}
          z(x) 
          \\ \hat u \\ \hat y \\ 0
               \end{bmatrix}
               + 
               \begin{bmatrix}
                    0 \\ 0 \\ I_m \\0
               \end{bmatrix}
               u, \\
               y &= \hat y,
 \end{align*}
 with $z(x)^\top=(z_1(x_1)^\top,\, z_2(x_2)^\top,\, \ldots,\, z_r(x_r)^\top)$,
 new dummy variables $\hat u, \hat y$ and setting 
 $X=\diag(X_1,X_2,\ldots,X_r)$ for $X \in \{E,J,R,B\} $.
 This coupling property of PH-DAE systems makes the port-Hamiltonian modelling framework 
  well suited for multiphysical applications.

Now, we consider external, time dependent inputs. To this end, we split the inputs and outputs into external (bar-notation) and internal ones (hat-notation), i.e.,  
 $B_i u_i$ is split into  $\bar B_i \bar u_i + \hat B_i \hat u_i$. Then, 
 the subsystem~\eqref{bcgjr_coupled.phdaes} reads 
 \begin{subequations}
 \label{bcgjr_coupled.phdaes2}
\begin{align}
     \ddt  (E_i x_i) & = (J_i-R_i) z_i(x_i) 
     + \hat B_i \hat u_i + \bar B_i \bar u_i,\\
     \hat y_i & = \hat B_i^{\top} z_i(x_i), \\
      \bar y_i & = \bar B_i^{\top} z_i(x_i). 
 \end{align}
   \end{subequations}
For the 
coupling relation (of the internal quantities) $\hat u+ C \hat y=0$ with a skew-symmetric matrix $C=-C^\top$ (which often arises in application), these systems can be written as a joint PH-DAE system in condensed form~\cite{gbjr21}:
 \begin{subequations}
 \begin{align}
     \ddt (Ex) & = (\tilde J-R) z(x) + \bar B \bar u, \\
     \bar y & = \bar B^\top z(x)
 \end{align}
  \end{subequations}
 with the condensed skew-symmetric matrix $\tilde J=J- \hat B C \hat B^\top$. Note that in this case all internal coupling modelled via the port-matrices $\hat B_i$ has now been transferred into the off-block diagonal elements of the skew-symmetric matrix $\tilde J$, i.e., $-\hat B C \hat B^\top$.
 
A~systems theoretic treatment of port-Hamiltonian systems goes back to {\sc Bernhard Maschke and Arjan van der Schaft} (see \cite{vdS04,JvdS14} for an overview), where nonlinear systems governed by ordinary differential equations are treated. For simplicity of presentation, we will (a) not follow the differential geometric path via Dirac structures, (b)  neglect a feed-through from input to output and (c) only consider finite dimensional systems, i.e., ordinary (ODEs) and differential-algebraic equations (DAEs), but no partial differential equations (PDEs). For simulation, the latter are usually transformed into ODEs and DAEs by spatial semi-discretization. For a differential geometric setting of PHS see~\cite{vdS06} and 
 an introduction into PH-PDEs see~\cite{JaSw12}. 
 
 The paper is organized as follows: In the next section we introduce PH-DAE systems which allow for a general nonlinear dissipative part. A PHS-DAE formulation of the MNA network equations is derived in Sect.~3, and for electromagnetic devices in Sect. 4. Section 5 discusses PHS formulation of coupled EM/circuit systems, which allow for monotolithic as well as weak coupling simulation approaches. Sect. 6  finishes with conclusions.
 
 \section{PH-DAE systems}
 When dealing with applications in electrical engineering, the concept of port-Hamiltonian modelling has to be generalized to coupled differential-algebraic equations, which (a) allow for a general nonlinear resistive part $r(z)$ instead of a quasilinear setting $R z$ as in the approach of~\cite{MeMo19} and (b) has only to be accretive on a subspace $\mathcal V \subset \mathbb R$ according to the constraints of the system.
 
A differential-algebraic equation of the form
\begin{equation}
\label{eq:def.phdae}
\begin{aligned}
\ddt E x(t)
    &= J z(x(t))-r(z(x(t))) +B u(t),
    \\
y(t) & =  B^{\top} z(x(t))
\end{aligned}
\end{equation}
is called a {\em port-Hamiltonian differential-algebraic equation} (PH-DAE)~\cite{gbjr21}, if the following holds:
\begin{itemize}
    \item $E\in\mathbb R^{n\times n}$, $\;J\in\mathbb R^{n\times n}$ and $B\in \mathbb R^{n\times m}$,
    $\;z,r:\mathbb R^{n}\to \mathbb R^{n}$.
    \item There exists a~subspace $\mathcal{V}\subset\mathbb R^n$ with the following properties:\\
    \begin{enumerate}
        \item[(i)] for all intervals $\interval\subset\mathbb R$ and functions $u:\interval\to\mathbb R^m$ such that \eqref{eq:def.phdae} has a~solution $x:\interval\to\mathbb R^n$, it holds $z(x(t))\in\mathcal{V}$ for all $t\in\interval$.
        \item[(ii)] $J$ is skew-symmetric on $\mathcal{V}$. That is, \quad
        $\displaystyle 
         v^\top J w=-w^\top J v %
        \;$   for all $\; v,w\in\mathcal{V}.$  
\item[(iii)] $r$ is accretive on $\mathcal{V}$. That is, 
        $\;v^\top r(v)\geq0 \;$ for all $\;v\in\mathcal{V}$.
    \end{enumerate}
\item There exists some function $H\in C^1(\mathbb R^n,\mathbb R)$ such that $\;\nabla H(x)=E^{\top}z(x) \;$ for all $\;x\in z^{-1}(\mathcal{V})$.
\end{itemize}
\begin{remark}
\begin{itemize}
\item[a)] The PH-DAE~\eqref{eq:def.phdae} system provides the usual energy balance
\begin{equation*}
\ddt H(x(t))  = - z(x(t))^\top r(z(x(t)))) + y(t)^\top u(t) \le  y(t)^\top u(t).
\label{eq:energbal}
\end{equation*}
\item[b)] PH-DAE subsystems now read
\begin{subequations}
\label{bcgjr_coupled.phdaes2}
\begin{align}
\ddt E_i  x_i(t) 
    =& J_i  z_i(x_i(t)) -r_i\bigl(z_i (x_i(t))\bigr) + B_i u_i (t), 
   \\
y_i(t) = &   B_i^\top z_i \bigl(x_i(t) \bigr)
 \end{align}
 \end{subequations}
instead of~\eqref{bcgjr_coupled.phdaes}, and if they are coupled by a skew-symmetric coupling relation  $\hat u+ C \hat y=0$ with a skew-symmetric matrix $C=-C^\top$ as before, they can be condensed into an overall PH-DAE system
\begin{subequations}
\label{eq:phdae.coupled.joint.network.condensed}
\begin{align}
\ddt E  x& =  \hat J z  - r +  \bar B \bar u, \\
\bar y & =  \bar B^{\top} z
\end{align}
\end{subequations}
with the skew-symmetric matrix $\hat J$ again given by  $\hat J = J - \hat B \hat C \hat B^{\top}$.
\end{itemize}
\end{remark}

\section{Electrical networks}

We consider the classical charge-/flux oriented MNA network equations \cite{Guenther1999,gbjr21}
 \begin{align*}
    \ddt\begin{bmatrix}
     0 & 0 & 0 &\Atop[C] & 0 \\
        0 & 0 & 0 & 0 & I \\
        0 & 0 & 0 & 0 & 0 \\
        0 & 0 & 0 & 0  & 0\\
        0 & 0 & 0 & 0  & 0
    \end{bmatrix}
    \begin{bmatrix}
     \npot{} \\ \jcurrv[L] \\ \jcurrv[V] \\ \charge[C] \\ \flux[L]   \end{bmatrix}
    = &
    \begin{bmatrix}
         0 & -\Atop[L] & -\Atop[V] & 0  & 0\\
         \Atop[L]^{\tr} & 0 & 0 & 0 & 0 \\
         \Atop[V]^{\tr} & 0 & 0 & 0 & 0 \\
         0 & 0 & 0 & 0 & 0\\
         0 & 0 & 0 & 0 & 0
    \end{bmatrix}
    \begin{bmatrix}
     \npot{} \\ \jcurrv[L] \\ \jcurrv[V] \\ \charge[C] \\ \flux[L]   \end{bmatrix}
    \\
   &  -
    \begin{bmatrix}
    \Atop[R]g(\Atop[R]^{\tr}\npot{})\\0\\ 
        0 \\ \charge[C]-\charge(\Atop[C]^{\tr} \npot{}) \\
        \flux[L]-  \flux(\jcurrv[L])
    \end{bmatrix}+ \begin{bmatrix}
    -\Atop[I] & 0  \\ 0 & 0 \\  0 & - I \\0 & 0 \\0 & 0 
    \end{bmatrix}
    \begin{bmatrix}
    \csrc({\tim}) \\ \vsrc({\tim})
    \end{bmatrix}\!
\end{align*}
with $e,\jmath_L$ and $\jmath_V$ denoting node potentials and currents through flux storing elements and voltages sources, $q_C$ and $\Phi_l$ charge and flux-storing elements, $i(t)$ and $v(t)$ independent current and voltage sources, the resistive currents $g$ and the incidence matrices $A_C,A_L,A_R,A_V,A_I$ for charge- and flux storing elements, resistive elements, voltage and current sources,    and seek a formulation as a PH-DAE system. For this, we need the following assumptions, which naturally occur in circuit simulation, see \cite{gbjr21}:
\begin{enumerate}
 \item[]
 \item[(a)] 
 \textbf{Soundness.} The circuit graph has at least one branch and is connected.
Furthermore, it 
contains neither $V$-loops nor $I$-cutsets. Equivalently, 
 $\Atop[V]$ and $\displaystyle(\Atop[C]\, \Atop[R]\, \Atop[L]\, \Atop[V] )^{\tr}$ have full column rank.
\item[(b)] \textbf{Passivity.} The functions $\charge$, $\flux$ and $\resi$ fulfill
\begin{enumerate}
\item[(i)] $\charge:\real^{n_C}\to\real^{n_C}$ and $\flux:\mathbbm{R}^{n_L}\to\mathbbm{R}^{n_L}$ are bijective, continuously differentiable, and their Jacobians 
\[
    \widetilde C(u_C):=\,\frac{\d \charge}{\d u_C}(u_C),
    \qquad
    \widetilde L(\jcurrv[L]):=\,\frac{\d \flux}{\d        
                \jcurrv[L]}(\jcurrv[L])
\]
are symmetric and positive definite for all $u_C\in \real^{n_C}$, $\jcurrv[L]\in \real^{n_L}$.
\item[(ii)] $\resi:\mathbbm{R}^{n_R}\to\mathbbm{R}^{n_R}$ is continuously differentiable, and its Jacobian has the property that $\frac{\d g}{\d u_R}(u_R)+ \frac{\d g}{\d u_R}(u_R)^{\tr}$ is positive definite for all $u_R\in\real^{n_R}$.
\end{enumerate}
\end{enumerate}
If $\charge:\real^{n_C}\to\real^{n_C}$ and $\flux:\mathbbm{R}^{n_L}\to\mathbbm{R}^{n_L}$ fulfill these assumptions, then there exist twice continuously differentiable and non-negative functions 
$V_C:\real^{n_C}\to\real$, $V_L:\real^{n_L}\to\real$ with the following property:
the gradients of $V_C$ and $V_L$ are, respectively, the inverse functions of $\charge$ and $\flux$. That is, 
\[
\begin{aligned}
&\forall \charge[C]\in\real^{n_C}: \nabla V_C(\charge[C])=\charge^{-1}(\charge[C]), 
 \,
&\forall \flux[L]\in\real^{n_L}: \nabla V_L(\flux[L])=\flux^{-1}(\flux[L]).
\end{aligned}\]
With this setting, the PH-DAE MNA network equations can now be derived as follows: we first eliminate the equation $\flux[L]-\flux(\jcurrv[L])=0$: $\jcurrv[L]=\flux^{-1}(\flux[L])$; secondly, we replace the equation $\charge[C]- \charge(\Atop[C]^{\tr} \npot)=0$ by $\Atop[C]^{\tr}\npot{}-\charge^{-1}(\charge[C])=0$. We end up with 
\begin{equation}
    \label{PH.DAE.MNA}
\begin{aligned}
    \ddt \underbrace{\begin{bmatrix}
     \Atop[C] & 0 & 0 & 0 \\
        0 & I & 0 & 0 \\
        0 & 0 & 0 & 0 \\
        0 & 0 & 0 & 0 \\
    \end{bmatrix}}_{\displaystyle  E:=}
    \underbrace{\begin{bmatrix}
    \charge[C] \\ \flux[L] \\ \npot{} \\ \jcurrv[V]     \end{bmatrix}}_{\displaystyle x:=}
    & =
    \underbrace{\begin{bmatrix}
         0 & -\Atop[L] & 0 & -\Atop[V] \\
         \Atop[L]^{\tr} & 0 & 0 & 0 \\
         0 & 0 & 0 & 0 \\
        \Atop[V]^{\tr} & 0 & 0 & 0 \\
    \end{bmatrix}}_{\displaystyle J:=}
    \underbrace{\begin{bmatrix}
    \npot{} \\  \flux^{-1}(\flux[L]) \\ \charge^{-1}(\charge[C]) \\  \jcurrv[V]
    \end{bmatrix}}_{\displaystyle z(x)}
    \\
    &
    -
    \underbrace{\begin{bmatrix}
    \Atop[R]g(\Atop[R]^{\tr}\npot{})\\0\\ 
         \Atop[C]^\tr \npot - \charge^{-1}(\charge[C]) \\0
    \end{bmatrix}}_{\displaystyle r(z(x)):=}+ \underbrace{\begin{bmatrix}
    -\Atop[I] & 0  \\ 0 & 0 \\ 0 & 0 \\ 0 & - I
    \end{bmatrix}}_{\displaystyle B:=}
    \underbrace{\begin{bmatrix}
    \csrc({\tim}) \\ \vsrc({\tim})
    \end{bmatrix}}_{\displaystyle u(t):=}\!,
\end{aligned}
\end{equation}
which is a PH-DAE of type~\eqref{eq:def.phdae} with subspace $\mathcal{V}$ and Hamiltonian $H(x)$ given by $H(x)= V_C(\charge[C]) + V_L(\flux[L]),  \mathcal{V}=\left\{\left.\begin{pmatrix}
        \npot{},\;  \jcurrv[L],\;  u_C,\;  \jcurrv[V]
    \end{pmatrix}^\tr \in\real^n\right|\Atop[C]^{\tr}e=u_C\right\}.$

\begin{remark}
\begin{itemize}
    \item[a)]
The PHS-DAE formulation shares the index properties of char\-ge/flux-oriented MNA network equations, if the assumption on soundness and passivity hold: the index is one if, and only if, it neither contains $LI$-cutsets nor $CV$-loops except for $C$-loops; otherwise it is two.
\item[b)] If $r$ subcircuits given as PH-DAE MNA network equations are coupled via voltage/current sources, the overall system can be written as a PH-DAE MNA of type~\eqref{PH.DAE.MNA}.
\end{itemize}
\end{remark}

\section{Electromagnetic devices}
In~\cite{Diab2022}, the Maxwell grid equations for an electromagnetic device have been developed as a linear PH-DAE system   provided that (a) the three-dimensional domain of the device is connected, bounded and surrounded by perfectly conducting material, (b) the permittivity $\epsilon$, the permeability $\mu$ are symmetric positive definite, and the conductivity $\sigma$ is symmetric positive semi-definite, and (c) finite integration technique \cite{Weiland96} has been used for the spatial discretization with orthogonal staggered cells: 
\begin{subequations}
\label{phs.model.em}
\begin{align}
    \begin{bmatrix} M_\mu  & 0 \\ 0 & M_\epsilon
    \end{bmatrix} \ddt \begin{bmatrix}
     \hat{h} \\  \hat{e}
    \end{bmatrix}& = \left( \begin{bmatrix}
    0 & -C \\ C^\top & 0
    \end{bmatrix}
    - \begin{bmatrix}
    0\; & 0 \\ 0\; & \;M_\sigma
    \end{bmatrix}\right) 
    \begin{bmatrix}
     \hat{h}  \\ \hat{e} 
    \end{bmatrix}
    + \begin{bmatrix} 
    0 \\ 
    X_S\color{black} 
    \end{bmatrix} \hat{u}_2, \\ 
    \hat y_2 & =  \begin{bmatrix} 
    0 \\ 
    X_S \color{black} 
    \end{bmatrix}^\top \begin{bmatrix} 
    \hat{h} \\ \hat{e} 
    \end{bmatrix} = 
    X_S^{\tr}  \color{black}\hat{e} .
\end{align}
\end{subequations}
Here $C$ denotes the discrete curl operator, the material matrices $M_\epsilon, M_\mu$ and $M_\sigma$ represent the discretized permittivity, permeability and conductivity distributions, $\hat{e}$ is vector of the electric mesh voltages $e$, $\hat{h}$ the vector of the magnetic mesh voltages $h$, and the (dual grid facet) source current $\hat{u}_2$ as input. 
This input is allocated at
positions $X_S$. In fact, $X_S$ maps the interior mesh links onto the exterior mesh nodes.
%
Furthermore, the respective electric mesh voltage $\hat y_2$ forms the output. \color{black} The Hamiltonian of the electromagnetic device is given by $H_1=\frac{1}{2} ( \tilde e^\top M_\epsilon \tilde e + \tilde h^\top M_\mu \tilde h ) $.

\section{Coupled EM/circuit system}
When coupling an electromagnetic device with an electric circuit, it remains only to define the inputs, outputs and the coupling equation. For the circuit, the electromagnetic device produces the current $\jmath_E$ flowing into the network, which is assembled at the respective nodes of the circuit via an incidence matrix $A_E$. Hence the circuit part reads (where we split inputs again in external inputs $\imath$, $v$, and internal ones):
\begin{subequations}
\label{phs.model.circuit}
\begin{align}
  \ddt \begin{bmatrix}
     \Atop[C] & 0 & 0 & 0 \\
        0 & I & 0 & 0 \\
        0 & 0 & 0 & 0 \\
        0 & 0 & 0 & 0 \\
        0 & 0 & 0 & 0 \\
    \end{bmatrix}
    \begin{bmatrix}
    \charge[C] \\ \flux[L] \\ \npot{} \\ \jcurrv[V] \\   \jcurrv[E]  \end{bmatrix}
    &=
    \begin{bmatrix}
         0 & -\Atop[L] & 0 & -\Atop[V] & - \Atop[E] \\
         \Atop[L]^{\tr} & 0 & 0 & 0 & 0\\
         0 & 0 & 0 & 0 & 0\\
        \Atop[V]^{\tr} & 0 & 0 & 0 & 0 \\
     \Atop[E]^{\tr} & 0 & 0 & 0 & 0 
    \end{bmatrix}
    \begin{bmatrix}
    \npot{} \\  \jcurrv[L]  \\ u_C \\  \jcurrv[V] \\ \jcurrv[E]
    \end{bmatrix}
    \\ \nonumber
    &
    -
    \begin{bmatrix}
    \Atop[R]g(\Atop[R]^{\tr}\npot{})\\0\\ 
         \Atop[C]^\tr \npot - u_C \\0 \\ 0 
    \end{bmatrix}+ 
        \begin{bmatrix}
    -\Atop[I] & 0 \\ 0 & 0  \\ 0 & 0  \\ 0 & - I  \\ 0 & 0 
    \end{bmatrix}
    \begin{bmatrix}
    \csrc({\tim}) \\ \vsrc({\tim}) 
    \end{bmatrix}
    +
    \begin{bmatrix}
    0 \\ 0   \\ 0  \\ 0  \\ 1 
    \end{bmatrix}
    \hat u_1, 
    \\ 
%
    \begin{bmatrix}
    {\bar y}_{1,1} \\ {\bar y}_{1,2} \\ \hat y_1 
    \end{bmatrix}
    & = \begin{bmatrix}
 -\Atop[I] & 0 & 0 \\ 0 & 0 & 0  \\ 0 & 0 & 0  \\ 0 & - I & 0 \\ 0 & 0 & 1 
    \end{bmatrix}^\top 
    \cdot \begin{bmatrix}
    \npot{} \\   \jcurrv[L] \\ u_C\\  \jcurrv[V] \\ \jcurrv[E]
    \end{bmatrix} = \begin{bmatrix}
    -\Atop[I]^{\tr} e \\-\jcurrv[V] \\  
    \jcurrv[E]
    \end{bmatrix}
\end{align}
\end{subequations}
with the Hamiltonian: $H_2= V_C(q_C)+V_L(\flux[L])$.

The coupling is as follows~\cite{Diab2022}: the input $\hat{u}_1$ (of the electric circuit) is given by the voltage drop at the electromagnetic device, which reads 
%
 
$\hat u_1 = - X_S^\top \tilde e = -\hat y_2$; 
on the other hand, the input $\hat u_2$ (of the magnetic device) is given by the current $\hat{u}_2= \jmath_E= \hat y_1$. 
Overall, we get the following skew-symmetric relation between inputs and outputs:
\label{phs.model.cc}
\begin{align}
\label{phs.model.cc}
0 & = 
\begin{bmatrix} \hat u_1 \\ \hat u_2  \end{bmatrix}
+
\begin{bmatrix}
0\;\; & I \\ - I\;\; & 0 \end{bmatrix}
\begin{bmatrix} \hat y_1 \\ \hat y_2  \end{bmatrix} .
\end{align}
\color{black}

As we have a system consisting of two PH-DAE systems~\eqref{phs.model.em} and~\eqref{phs.model.circuit} with a skew-symmetric linear coupling condition~\eqref{phs.model.cc}, the overall system can be written as a condensed PH-DAE system~\eqref{eq:phdae.coupled.joint.network.condensed} with Hamiltonian $H=H_1+H_2$ and enlarged matrices as above.
\section{Simulation Strategies}
Generally, for simulating the coupled EM/circuit system numerically, two approaches are feasible:
\begin{itemize}
    \item[a)] {\em Monolithic approach.} The condensed system~\eqref{eq:phdae.coupled.joint.network.condensed} can be solved by any integration scheme suitable for index-1 and index-2 systems, depending on the index. To preserve the dissipation inequality also on a discrete level, collocation schemes~\cite{MeMo19} and discrete gradient schemes tracing back to~\cite{Gonz96}  are the methods-of choice. This strategy is also referred to as strong coupling. 
    \item[b)] 
    %
    %
    {\em Monolithic multirate approach.} In fact, we are facing models, where the subsystems can have widely separated time scales. This can create so-called multirate potential, where it is beneficial to employ schemes, which use inherent step sizes for each subsystem. In this way, each subsystem can be sampled on its time scale. See e.g. \cite{Guenther2016,Bartel2022}.
    \color{black}
    \item[c)] {\em Weak coupling.}  Since the coupling equations is merely the one-to-one identification of output and input, we
    can insert this. Furthermore, omitting outputs due to external sources, we have
    \begin{subequations}
    \label{circuit.em.coupled2}
    \begin{align}
    \ddt 
    &
    \begin{bmatrix}
     \Atop[C] & 0 & 0 & 0 \\
        0 & I & 0 & 0 \\
        0 & 0 & 0 & 0 \\
        0 & 0 & 0 & 0 \\
        0 & 0 & 0 & 0 \\
    \end{bmatrix}
    \begin{bmatrix}
    \charge[C] \\ \flux[L] \\ \npot{} \\ \jcurrv[V] \\   \jcurrv[E]  \end{bmatrix}
    =
    \begin{bmatrix}
         0 & -\Atop[L] & 0 & -\Atop[V] & - \Atop[E] \\
         \Atop[L]^{\tr} & 0 & 0 & 0 & 0\\
         0 & 0 & 0 & 0 & 0\\
        \Atop[V]^{\tr} & 0 & 0 & 0 & 0 \\
     \Atop[E]^{\tr} & 0 & 0 & 0 & 0 
    \end{bmatrix}
    \begin{bmatrix}
    \npot{} \\  \jcurrv[L]  \\ u_C \\  \jcurrv[V] \\ \jcurrv[E]
    \end{bmatrix}
    \\ \nonumber
    & \hspace*{2cm}
    -
    \begin{bmatrix}
    \Atop[R]g(\Atop[R]^{\tr}\npot{})\\0\\ 
         \Atop[C]^\tr \npot - u_C \\0 \\ 0 
    \end{bmatrix}+ 
        \begin{bmatrix}
    -\Atop[I] & 0 \\ 0 & 0  \\ 0 & 0  \\ 0 & - I  \\ 0 & 0 
    \end{bmatrix}
    \begin{bmatrix}
    \csrc({\tim}) \\ \vsrc({\tim}) 
    \end{bmatrix}
    -
    \begin{bmatrix}
    0 \\ 0   \\ 0  \\ 0  \\ 1 
    \end{bmatrix}
    \hat y_2,
    \\ 
%
   & \hspace*{3cm} \hat y_1 
    = 
    \jcurrv[E]
\end{align}      
    \end{subequations}
    and
    \begin{subequations}
\label{phs.model.em2}
\begin{align}
    \begin{bmatrix} M_\mu  & 0 \\ 0 & M_\epsilon
    \end{bmatrix} \ddt \begin{bmatrix}
     \tilde h \\  \tilde e
    \end{bmatrix} = &\left( \begin{bmatrix}
    0 & -C \\ C^\top & 0
    \end{bmatrix}
    - \begin{bmatrix}
    0 & 0 \\ 0 & M_\sigma
    \end{bmatrix}\right) 
    \begin{bmatrix}
     \tilde h \\ \tilde e
    \end{bmatrix}
    + 
    \begin{bmatrix} 0 \\ X_S\end{bmatrix} \hat{y}_1
    \color{black}
    \\ %
    \hat{y}_2 = &     
             X_S \tilde{e}. \color{black}
\end{align}
\end{subequations}
Here dynamic iteration schemes~\cite{ArGu01} are the methods-of choice, as due to the ODE-DAE coupling no stability constraints occur~\cite{BBGS13}. In addition, each step of a Jacobi or Gau\ss-Seidel iteration scheme defines a PH-DAE system by its 
own~\cite{gbjr21}. 

Operator splitting approaches are not generally feasible for differential-algebraic equations, which can easily be seen for the linear PH-DAE~\eqref{eq.lin.ph-dae} with $z(x)=x$ and $B=0$. A Lie-Trotter splitting approach may read
\begin{align*}
\ddt (Ex) & = Jx, \quad x(0)=x_0, \\
\ddt (Ew) & = - Rw. \quad w(0)=x(T),
\end{align*}
allowing for using a symplectic integrator for the first step, and a dissipative scheme for the second one. However, the matrix pencil $\{E,J\}$ or $\{E,R\}$ may be singular and thus not define a unique solution for the respective subproblem, even if   the matrix pencil $\{E,J-R\}$ of the overall system is regular.  Even if this does not happen,
the first problem, for example, may not allow for a unique solution for arbitrary choices of consistent initial values. For
\begin{align*}
    E= \diag(1,0,1), \quad J=\begin{bmatrix}
         0 & -1 & 0 \\ 1 & 0 & 0 \\ 0 & 0 & 0
    \end{bmatrix}, \quad
    R=\diag(0,1,1), \quad x_0=\begin{bmatrix}
         1 \\ -1 \\ 0
    \end{bmatrix},
\end{align*}
all matrix pencils $\{E,J-R\}$, $\{E,J\}$ and $\{E,R\}$ are regular, but the first step yields $x_1\equiv 0 \neq 1$.

One may overcome the problem by rewriting the DAE in terms of an underlying ODE and subsequent algebraic variables given by explicit evaluations.
For network equations a branch oriented
loop-cutset approach is an option for defining such a PH-DAE system, see~\cite{Diab2022}. Another way 
to avoid the problems above is to follow an 
operator splitting based approach for dynamic iteration. In the latter case, no stability problems occur and a monotone convergence can be obtained~\cite{BGJR22}. 
\end{itemize}
 
\section{Conclusions} PHS provide a modelling framework which preserves essential physical properties. It is especially suited for multiphysical applications, as the proper coupling of port-Hamiltonian subsystems yields an overall PHS. In electrical engineering, we have shown that electrical networks and electromagnetic devices can be written as PHS, and coupled EM/circuit system yield coupled PHS with a skew-symmetric coupling, which can be rewritten as an overall PHS. For simulation, a monolithic approach is suitable for the former, and weak coupling methods for the latter. There are still many unresolved questions, such as how to adequately integrate distributed ports into the PHS modeling.

\begin{acknowledgement}
Michael Günther is  indebted to the funding given by the European Union’s Horizon 2020 research and innovation programme under the Marie Sklodowska-Curie Grant Agreement No. 765374, ROMSOC.
\end{acknowledgement}



\end{document}